\documentclass[12pt, reqno]{amsart}
\usepackage{amsmath}
\usepackage{amssymb}
\usepackage{amsfonts}
\usepackage{mathrsfs}
\usepackage{bbding}
\usepackage{amssymb,amsmath,amscd,graphicx,latexsym}
\newtheorem{theorem}{Theorem}[section]

\newtheorem{proposition}[theorem]{Proposition}
\newtheorem{lemma}[theorem]{Lemma}

\renewcommand{\a}{\alpha}
\renewcommand{\b}{\beta}
\renewcommand{\d}{\delta}

\newcommand{\f}{\frac}
\newcommand{\g}{\gamma}

\newcommand{\bea}{\begin{eqnarray}}
\newcommand{\eea}{\end{eqnarray}}
\newcommand{\bna}{\begin{eqnarray*}}
\newcommand{\ena}{\end{eqnarray*}}
\renewcommand{\o}{\omega}
\renewcommand{\O}{\Omega}

\newcommand{\s}{\sigma}
\renewcommand{\th}{\theta}

\newcommand{\ve}{\varepsilon}

\begin{document}

\title[Multiplicative functions resembling the M\"obius funciton]  
{Multiplicative functions resembling \\ the M\"obius funciton}  
\author{Qingyang Liu}
\address{School of Mathematics
\\
Renmin University of China 
\\
Beijing 100872
\\
China}

\address{Institute of Mathematics \\ 
AMSS\\ 
The Chinese Academy of Sciences\\ 
Beijing 100080
\\
China}
\email{qingyangliu@amss.ac.cn} 

\begin{abstract}
A multiplicative function $f$ is said to be resembling the M\"{o}bius function 
if $f$ is supported on the square-free integers, and $f(p)=\pm 1$ for each prime $p$. 
We prove $O$- and $\O$-results for the 
summatory function $\sum_{n\leq x} f(n)$ for a class of these $f$ studied by Aymone, 
and the point is that these $O$-results demonstrate cancellations better than the square-root saving. It is proved 
in particular that the summatory function is $O(x^{1/3+\ve})$ under the Riemann Hypothesis. 
On the other hand it is proved to be $\Omega(x^{1/4})$ unconditionally. It is interesting 
to compare these with the corresponding results for the M\"{o}bius function.  
\end{abstract}

\date{\today}


\maketitle

\addtocounter{footnote}{1}

\section{Introduction}
\setcounter{equation}{0}
 
The M\"{o}bius function $\mu$ is an arithmetic multiplicative function supported 
on square-free integers with $\mu(p)=-1$ for each prime $p$. Thus $\mu(1)=1$, 
$\mu(n)=0$ if $n$ is divided by a square of prime, and $\mu(n)=(-1)^k$ if 
$n$ is a product of $k$ distinct prime factors. The role of $\mu$ is central in the theory of prime 
numbers, and an evidence is that the Riemann Hypothesis is equivalent to the estimate 
$M(x)\ll x^{1/2+\ve}$ for arbitrary $\ve>0$, where $M(x)$ denotes the summatory function of $\mu$, that is
\bna 
M(x)=\sum_{n\leq x} \mu(n). 
\ena  
It is also known, unconditionally, that $M(x)=\Omega(\sqrt{x})$ 
with the $\O$-symbol being defined as the negation of the $o$-symbol. 
All these are classical, and the reader is referred to Titchmarsh \cite[Chap. 14]{Tit} for details. 
The best-known upper bound estimate for $M(x)$ is obtained in Soundararajan \cite{Sou}. 

A random model of the M\"{o}bius function was proposed by Wintner \cite{Win}, and 
studied by many authors from various perspectives. An arithmetic function $f$ is said to be a {\it random multiplicative function} 
if 
\begin{itemize}
\item[(i)] $f$ is multiplicative, and supported on the square-free integers;  
\item[(ii)] $(f(p))_p$ is an independent sequence of random variables with distribution 
$$
\mathbb{P}(f(p)=1) = \mathbb{P}(f(p)=-1) = \f12, 
$$ 
where $p$ runs over the set of primes. 
\end{itemize}  
Denote by $M_f(x)$ the summatory function 
\bna 
M_f(x) =\sum_{n\leq x} f(n). 
\ena  
Then a theorem of Wintner states that $M_f(x)\ll x^{1/2+\ve}$ almost surely, and also $M_f(x)=\O (x^{1/2-\ve})$ almost 
surely. For developments concerning the above $\ll$- and $\O$-results of $M_f$, 
see Hal\'{a}sz \cite{Hal}, 
Lau-Tenenbaum-Wu \cite{LauTenWu},  
and Harper \cite{Har} and the references therein.      

An arithmetic function $f$ is said to be {\it resembling the M\"{o}bius function} 
if 
\begin{itemize}
\item[(i)] $f$ is multiplicative, and supported on the square-free integers, 
\item[(ii)] $f(p)=\pm 1$  for each prime $p$.  
\end{itemize}  
It is therefore an interesting problem whether there is an $f$ resembling the M\"{o}bius function, such that 
its summatory function satisfies $M_f(x)=o(\sqrt{x})$. Aymone \cite{Aym} has studied this problem and provided   
a class of examples of these $f$. Let $\chi$ be a real non-principal Dirichlet character modulo $q\geq 3$, and define 
\bea\label{def/gchi} 
g_\chi (p) = 
\left\{
\begin{array}{lll}
\chi(p), & \text{ if } p\nmid q;  \\
1, & \text{ if } p|q. 
\end{array}
\right. 
\eea 
This $g_\chi$ extends to 
a completely multiplicative function $g_\chi: \mathbb{N} \to \{-1, 1\}$.    
The example discovered in \cite{Aym} is $f=\mu^2 g_\chi$, where clearly $f$ is resembling the M\"{o}bius function. 
It is proved in \cite{Aym} that, under the condition 
\bea\label{con|chi} 
\sum_{p\leq x} |1-f(p)\chi(p)|\ll \sqrt{x} \exp(-c\sqrt{\log x}), 
\eea 
one has 
\bea\label{x|1/4} 
M_f(x) \ll \sqrt{x} \exp(-c(\log x)^{1/4})    
\eea 
where here and throughout $c$ stands for positive constants not necessarily the 
same at each occurrence; 
while under the Riemann Hypothesis for the zeta-function, 
\bea\label{x|2/5}
M_f(x)\ll x^{2/5+\ve}
\eea  
for arbitrary $\ve>0$. 

\medskip 

The purpose of this paper is to improve on these, and also prove an $\O$-result. Our main results 
are stated in the following two theorems. 

\begin{theorem}\label{thm1} 
Let $\chi$ be a real non-principal Dirichlet character modulo $q\geq 3$, and 
$g_\chi$ the completely multiplicative function extended from $\chi$ as in \eqref{def/gchi}. 
Let $f=\mu^2 g_\chi$. 
Then 
\begin{itemize} 
\item[(i)] there exists a positive constant $c$ such that   
\bna
M_f(x)\ll x^{1/2} \exp\bigg(-c \f{(\log x)^{3/5}}{(\log\log x)^{1/5}}\bigg); 
\ena  
\item[(ii)] under the Riemann Hypothesis for the zeta-function, 
\bna
M_f(x)\ll x^{1/3+\ve}. 
\ena  
\end{itemize} 
\end{theorem}

Theorem~\ref{thm1} improves \eqref{x|1/4} without applying the additional condition \eqref{con|chi}, and improves 
\eqref{x|2/5} under the same condition. Theorem~\ref{thm1}  can be compared with 
the $\O$-result below. 

\begin{theorem}\label{thm3} 
Let $\chi$ be a real non-principal Dirichlet character modulo $q\geq 3$, and 
$g_\chi$ the completely multiplicative function extended from $\chi$ as in \eqref{def/gchi}. 
Let $f=\mu^2 g_\chi$. 
Then 
\bna
M_f(x)=\O(x^{1/4}). 
\ena  
\end{theorem}
Note that Theorem~\ref{thm3} is unconditional. 
The same result has been proved by Aymone \cite{Aym} under the Riemann Hypothesis for $L(s,\chi)$. 
Unconditionally Klurman et al \cite{Klur} have proved that $M_f(x)=\O(x^{1/4-\ve})$ for arbitrary $\ve>0$. 

\medskip 
 
We use standard notations in number theory. The Riemann Hypothesis for the zeta-function means that 
all the non-trivial zeros of $\zeta(s)$ lie on the critical line $\s=\f12$. The conductor $q$ of $\chi$ is 
considered as fixed, and therefore some $\ll$ or $O$-constants may depends on $q$, 
but we do not make these 
dependences explicit. The expression $f(x)=\O(g(x))$ means the negation of the estimate $f(x)=o(g(x))$, 
that is the inequality $|f(x)|\geq c g(x)$ is satisfied for some arbitrarily large values of $x$. 
The letter $c$ stands for positive constants not necessarily the 
same at each occurrence, and $s=\s+it$ denotes a complex variable.  
The letter $\ve>0$ stands for a real number arbitrarily small, not necessarily the 
same at each occurrence. 
 
\section{Preparations}    
\setcounter{equation}{0}

To prove Theorems~\ref{thm1} and \ref{thm3}, the first step is to compute 
the generating function for $f$. 
Recall that $f=\mu^2 g_\chi$ and $g_\chi$ 
is the completely multiplicative function extended from $\chi$ as in \eqref{def/gchi}.   
Write $s=\s+it$ as usual. Then in the half-plane $\s>1$ we have    
\bna 
\sum_{n=1}^\infty \f{f(n)}{n^s} 
&=& \prod_p \bigg(1+ \f{f(p)}{p^s}\bigg) 
= \prod_p \bigg(1+ \f{g_\chi(p)}{p^s}\bigg) \\ 
&=& \prod_{p|q} \bigg(1+ \f{1}{p^s}\bigg) \prod_{p\nmid q} \bigg(1+ \f{\chi(p)}{p^s}\bigg). 
\ena 
Since $\chi$ is real, we have $\chi(p)^2=1$ for $p\nmid q$, and hence  
\bna 
\bigg(1+ \f{\chi(p)}{p^s}\bigg)
=\bigg(1- \f{\chi(p)}{p^s}\bigg)^{-1} \bigg(1- \f{1}{p^{2s}}\bigg). 
\ena 
It follows that 
\bea\label{sumfn=} 
\sum_{n=1}^\infty \f{f(n)}{n^s} 
= \f{L(s, \chi) P(s)}{\zeta(2s)} 
\eea 
where 
\bea\label{Ps=} 
P(s)= \prod_{p|q} \bigg(1- \f{1}{p^s}\bigg)^{-1}. 
\eea 
The function $P$ above has no zero, and has infinitely many poles at 
\bea\label{P/pol} 
\rho = i \f{2\pi j}{\log p}, \quad p|q, \  j=0, \pm 1, \pm2, \ldots 
\eea 
which are lying on the imaginary axis. The pole with $j=0$ is of order $\o(q)$ where $\o(q)$ denotes the 
number of distinct prime divisors of $q$, while other poles are simple. 
The formula \eqref{sumfn=} holds in the half-plane $\s>1$, while the function on the right-hand side of \eqref{sumfn=} 
is meaningful in the whole complex plane by the functional equations 
of $\zeta(s)$ and $L(s,\chi)$. 

This section is devoted to lemmas that are needed for Theorem~\ref{thm1}. 

\begin{lemma}\label{lem/h} 
Let $P$ be as in \eqref{Ps=}, and $h(n)$ defined as the coefficients in the Dirichlet series 
expression  
\bea\label{def/h}
\f{P(s)}{\zeta(2s)}=\sum_{n=1}^\infty \f{h(n)}{n^{s}} 
\eea
which holds in the half-plane $\s>\f12$. Denote by $M_h(x)$ the summatory function $\sum_{n\leq x} h(n)$. 
Then there exists a positive constant $c$ such that   
\bna
M_h(x)\ll x^{1/2} \exp\bigg(-c \f{(\log x)^{3/5}}{(\log\log x)^{1/5}}\bigg). 
\ena  
\end{lemma}

\begin{proof} By definition \eqref{def/h},  
\bna 
h(n) = \sum_{dm^2=n\atop p|d\Rightarrow p|q} \mu(m) 
\ena 
where the condtition ``$p|d\Rightarrow p|q$'' means that all the prime divisors of $d$ are divisors of $q$. 
It follows that 
\bna 
M_h(x)= \sum_{d\leq x\atop p|d\Rightarrow p|q} \sum_{m^2\leq x/d} \mu(m). 
\ena 
Now we invoke the estimate (see for example \cite[Theorem 12.7]{Ivi}) 
\bea\label{VinKor} 
\sum_{m\leq y} \mu(m) \ll y \exp\bigg(-c \f{(\log y)^{3/5}}{(\log\log y)^{1/5}}\bigg)
\eea 
which follows from the Vinogradov-Korobov zero-free region of the Riemann zeta-function. 
Writing $\Phi(y)$ for the function on the right-hand side of \eqref{VinKor}, we see that $\Phi(y)$ 
is increasing for $y\geq 100$, say.  Thus, for $\sqrt{x/d}\geq 100$, 
\bna
\sum_{m^2\leq x/d} \mu(m)  
\ll \Phi\bigg(\sqrt{\f{x}{d}}\bigg) \ll \Phi(\sqrt{x}). 
\ena 
The above sum is bounded  for $\sqrt{x/d}\leq 100$, and therefore 
\bea\label{smPhi} 
M_h(x)\ll \Phi(\sqrt{x}) \sum_{d\leq x\atop p|d\Rightarrow p|q} 1. 
\eea 

We are going to prove that 
\bea\label{sum/dypd} 
\sum_{d\leq x \atop p|d\Rightarrow p|q} 1 \leq (\log x)^{\o(q)} 
\eea 
where $\o(q)$ is the number of distinct prime divisors of $q$. 
We let $\o(q)=k$ and let $q=p_1^{\a_1}\cdots p_k^{\a_k}$ be the canonical decomposition of $q$. 
The condition ``$p|d\Rightarrow p|q$'' means that $d$ must be of the form 
$d=p_1^{\b_1}\cdots p_k^{\b_k}$ 
where $\b_1, \ldots, \b_k$ are nonnegative integers. Note that $\b_j$ may exceed $\a_j$, but 
they must be bounded from above as 
\bna 
\b_j\leq \f{\log x}{\log p_j} \leq \log x, 
\ena 
since $p_j^{\b_j}\leq x$ for each $j$. This proves \eqref{sum/dypd}. 
The assertion of the lemma is a consequence of \eqref{smPhi} and \eqref{sum/dypd}. 
\end{proof} 

\begin{lemma}\label{lem/|h|} 
Let $h(n)$ be as in Lemma~\ref{lem/h}. Then 
\bna 
\sum_{n\leq x} |h(n)| \ll x^{1/2} (\log x)^{\o(q)}, 
\ena 
where $\o(q)$ denotes the number of different prime divisors of $q$. 
\end{lemma}

\begin{proof} 
By \eqref{def/h},  
\bea\label{|h|} 
|h(n)|=\bigg|\sum_{dm^2=n\atop p|d\Rightarrow p|q} \mu(m)\bigg|
\leq 
\sum_{dm^2=n\atop p|d\Rightarrow p|q} |\mu(m)|. 
\eea 
Taking summation over $n$, we have 
\bna 
\sum_{n\leq x} |h(n)| \ll 
\sum_{n\leq x} \sum_{dm^2=n\atop p|d\Rightarrow p|q} |\mu(m)| 
 \ll \sum_{m\leq \sqrt{x}} \sum_{d\leq x/m^2 \atop p|d\Rightarrow p|q} 1.  
\ena 
The assertion of the lemma now follows from this and \eqref{sum/dypd}.  \end{proof}

\section{Proof of Theorem~\ref{thm1}(i)}
\setcounter{equation}{0}
     
\begin{proof}[Proof of Theorem~\ref{thm1}(i)] 
We start from \eqref{sumfn=}  and \eqref{def/h} to get 
\bna 
f(n) = \sum_{dm=n} \chi(d) h(m), 
\ena 
and hence 
\bea\label{Mfx=SSS} 
M_f(x)
&=&\sum_{n\leq x} f(n) 
= \sum_{dm\leq x} \chi(d) h(m) \nonumber\\ 
&=& \bigg(\sum_{dm\leq x\atop d\leq D} + \sum_{dm\leq x\atop m\leq M} - \sum_{d\leq D \atop m\leq M} \bigg) 
\chi(d) h(m)\nonumber\\
&=& S_1 + S_2 - S_3, 
\eea 
say. Here $D$ and $M$ are parameters satisfying $DM=x$ but to be decided later.  
We start from 
\bna
S_1
\ll \sum_{d\leq D} \bigg|\sum_{m\leq x/d} h(m) \bigg|, 
\ena 
and then apply Lemma~\ref{lem/h} to the inner sum. 
Note that the function 
\bna 
\exp\bigg(-c \f{(\log x)^{3/5}}{(\log\log x)^{1/5}}\bigg)
\ena   
is decreasing. By Lemma~\ref{lem/h}, 
\bea\label{S1=} 
S_1 
&\ll& \sum_{d\leq D} \bigg(\f{x}{d}\bigg)^{1/2} 
\exp\bigg(-c \f{(\log \f{x}{D})^{3/5}}{(\log\log \f{x}{D})^{1/5}}\bigg) \nonumber\\  
&\ll& D^{1/2} x^{1/2} 
\exp\bigg(-c_1 \f{(\log \f{x}{D})^{3/5}}{(\log\log \f{x}{D})^{1/5}}\bigg)
\eea 
for some positive constant $c_1$. 
Now we invoke the trivial bound $|\sum_{d\leq x}\chi(d)|\leq q$ as well as  Lemma~\ref{lem/|h|}, to get 
\bea\label{S2=} 
S_2
= \sum_{m\leq M} h(m) \sum_{d\leq x/m} \chi(d)   
\ll \sum_{m\leq M} |h(m)| 
\ll M^{1/2} \log^{\o(q)} x, 
\eea 
where we recall that in this paper $q$ is treated as a constant. Furthermore, 
by Lemma~\ref{lem/h}, 
\bea\label{S3=} 
S_3
= \bigg(\sum_{m\leq M} h(m)\bigg)\bigg(\sum_{d\leq D} \chi(d)\bigg)
\ll \bigg|\sum_{m\leq M}h(m) \bigg|
\ll M^{1/2}. 
\eea 
It turns out that the optimal choice is  
\bna 
M= x 
\exp\bigg(-\f{c_1}{2} \f{(\log x)^{3/5}}{(\log\log x)^{1/5}}\bigg)
\ena 
and $D=x/M$. This proves Theorem~\ref{thm1}(i).  
\end{proof} 

\section{Further preparations} 
\setcounter{equation}{0}

This section is devoted to lemmas that are needed for Theorem~\ref{thm1}(ii).  
The first of these is  Perron's formula in the following form, which is a modification of 
Titchmarsh \cite[Lemma~3.12]{Tit}. 

\begin{lemma}\label{lem1} 
Let $F(s)=\sum a_n n^{-s}$ be the Dirichlet series expression for $F$ in the half-plane $\s>1$, 
where $a_n\ll \psi(n)$, $\psi$ being non-decreasing, and 
\bea\label{|an|<<} 
\sum_{n=1}^\infty \f{|a_n|}{n^\s} \ll \f{1}{\s-1}
\eea  
as $\s\to 1^+$. If $b>0, b+\s>1, x$ is half an odd integer, then 
\bea\label{sum|an} 
\sum_{n\leq x} \f{a_n}{n^s} 
&=& \f{1}{2\pi i} \int_{b-iT}^{b+iT} F(s+w) \f{x^w}{w}dw +O\bigg(\f{x^b}{T(\s+b-1)}\bigg) \nonumber\\ 
&& +O\bigg(\f{\psi(2x)x^{1-\s} \log x}{T}\bigg). 
\eea   
\end{lemma} 

\begin{lemma}\label{lem20} 
Assume the Riemann Hypothesis for the zeta-function. Let $m\geq 0$ be a fixed integer. 
Then in the half-plane $\s\geq \f12+\d$ we have 
\bna 
\zeta^{(m)}(s) \ll (|t|+2)^\ve
\ena 
and 
\bna 
\f{d^m}{ds^m}\bigg(\f{1}{\zeta(s)}\bigg)\ll (|t|+2)^\ve 
\ena 
where $\d>0$ and $\ve>0$ are arbitrarily small. We may require that $\d\leq \ve$. 
\end{lemma}

\begin{proof} 
This follows from \cite[Theorem~14.14(B)]{Tit} and \cite[Theorem~14.16]{Tit}. 
\end{proof} 

\begin{lemma}\label{lem21} 
Unconditionally we have, uniformly for $\f12\leq \s\leq 2$,  
\bna 
\int_{-T}^T \f{|L(\s+it,\chi)|}{|\s+it|} dt \ll \log^2 T. 
\ena 
\end{lemma}

\begin{proof} By Cauchy's inequality, 
\bna 
\int_{-T}^T \f{|L(\s+it,\chi)|}{|\s+it|} dt 
\ll \bigg(\int_{-T}^T \f{1}{|\s+it|} dt\bigg)^{1/2}\bigg(\int_{-T}^T \f{|L(\s+it,\chi)|^2}{|\s+it|} dt \bigg)^{1/2}. 
\ena 
The first integral in braces is $\ll \log T$. To bound the second integral in braces, we note that, 
uniformly for $\f12\leq \s\leq 2$, 
\bna 
\int_{-T}^T |L(\s+it,\chi)|^2 dt \ll T\log T. 
\ena 
It follows by partial integration that 
\bna 
\int_{-T}^T \f{|L(\s+it,\chi)|^2}{|\s+it|} dt 
\ll \log^2 T. 
\ena 
This proves the lemma. \end{proof}

The above lemmas will be applied to estimate the tail in \eqref{def/h}, that is 
\bea\label{def/ZM}
Z_M(s)
= \f{P(s)}{\zeta(2s)} -\sum_{m\leq M} \f{h(m)}{m^{s}} 
\eea
with $M\geq 1$ being a parameter, in the half-plane $\s\geq \f12+\d$. We observe that 
\bea\label{ZM|>1}
Z_M(s)
=\sum_{m>M} \f{h(m)}{m^{s}}
\eea
in the half-plane $\s\geq 1+\d$. In the above $\d>0$ is arbitrarily small.  

\begin{proposition}\label{prop/ZM}  
Assume the Riemann Hypothesis for the zeta-function. Let $Z_M(s)$ be as in \eqref{def/ZM}.  
Then, in the half-plane $\s\geq \f12+\d$ with $\d>0$ arbitrarily small, 
\bna
Z_M(s)\ll M^{1/4-\s} M^\ve (|t|+2)^{\ve}. 
\ena
\end{proposition}

\begin{proof} 
By definition \eqref{Ps=},   
\bna 
P(s)=\sum_{d\atop p|d\Rightarrow p|q} \f{1}{d^s} 
\ena 
where the condtition ``$p|d\Rightarrow p|q$'' means that all the prime divisors of $d$ are divisors of $q$. 
Hence by \eqref{def/h},  
\bea\label{|h|} 
|h(n)|=\bigg|\sum_{dm^2=n\atop p|d\Rightarrow p|q} \mu(m)\bigg|
\leq 
\sum_{dm^2=n\atop p|d\Rightarrow p|q} |\mu(m)|. 
\eea 
If we relax the constraint on $d$, we have
$$
|h(n)|\leq \sum_{m^2|n}|\mu(m)| \leq 2^{\o(n)}, 
$$
where $\o$ is the number of distinct prime divisors of $n$. 
A rough bound for the above is $|h(n)|\ll n^\ve$.   
It follows from \eqref{|h|} that 
\bea\label{sum/|h|} 
\sum_{n=1}^\infty \f{|h(n)|}{n^{\s}} 
\leq 
P(\s) \sum_{n=1}^\infty \f{|\mu(n)|}{n^{2\s}} 
\ll \sum_{n=1}^\infty \f{|\mu(n)|}{n^{\s}} 
=\f{\zeta(\s)}{\zeta(2\s)}
\ll \f{1}{\s-1}
\eea 
as $\s\to 1$ from the right.  

We apply Perron's formula in the form of Lemma~\ref{lem1}, in which we  
take $\psi(n)=n^\ve$ and $b= 1 +\f{1}{\log M}$. Then 
\eqref{sum|an} becomes, for $M$ being half an odd integer,  
\bea\label{Mf/int/V} 
\sum_{m\leq M} \f{h(m)}{m^{s}}
= \f{1}{2\pi i} \int_{b-iV}^{b+iV} \f{P(s+w)}{\zeta(2s+2w)} \f{M^w}{w} dw +O\bigg(\f{M^{1+\ve}}{V}\bigg). 
\eea 
Write $w=u+iv$ where $u$ and $v$ are the real and imaginary parts of $w$ respectively. 
We shift the contour in \eqref{Mf/int/V} to the left till the vertical line $u=a$ where
\bea\label{a/ZM} 
a=\f14-\s +\d 
\eea 
with $\d>0$ arbitrarily small, passing the simple pole of the integrand at $w=0$, and also 
poles at $w=\rho-s$ with $\rho$ as in \eqref{P/pol} satisfying $|v|\leq V$. 
Hence the integral in \eqref{Mf/int/V} can be written as 
\bea\label{V/iii} 
\int_{b-iV}^{b+iV} 
&=& - \int_{C_1} + \int_{C_2} + \int_{a-iV}^{a+iV}  + R
\eea 
where $C_1$ means the segment $\{w=u-iV: a\leq u\leq b\}$, $C_2$ means $\{w=u+iV: a\leq u\leq b\}$,  
and $R$ denotes the summation of all the residues arising from the poles mentioned above.  

We want to estimate the three integrals in \eqref{V/iii}, and to this end we need to estimate 
their integrands in the half-plane $u\geq a$. Note that we are assuming that $\s\geq \f12+\d$. 
Thus, in the half-plane $u\geq a$ under consideration, we have 
$\Re(w+s)=u+\s\geq a+\s=\f14+\d$, and therefore  
\bna 
P(w+s)
\ll \prod_{p|q} \bigg(1-\f{1}{p^{1/4+\d}}\bigg)^{-1} 
\ll \prod_{p|q} \f{p^{\d}}{p^{\d}-1}\ll q^{\d}.    
\ena 
The definition \eqref{a/ZM} of $a$ gurantees that, in the half-plane $u\geq a$,  
we have $\Re (2s+2w)=2\s+2u\geq \f12+\d$, and Lemma \ref{lem20} applies. 
Hence the integrals on $C_1$ and $C_2$ 
can be estimated as  
\bna 
\bigg(-\int_{C_1} + \int_{C_2}\bigg) \f{P(s+w)}{\zeta(2s+2w)} \f{M^w}{w} dw 
\ll \int_a^b (|t|+V+2)^\ve \f{M^u}{V} du 
\ll M V^{\ve-1} (|t|+2)^\ve.   
\ena 
On the vertical segment, we have 
\bna 
\int_{a-iV}^{a+iV} \f{P(s+w)}{\zeta(2s+2w)} \f{M^w}{w} dw  
&\ll& M^a    
\int_{-V}^V \f{(|t|+|v|+2)^\ve}{\sqrt{a^2+v^2}} dv \\ 
&\ll& M^{1/4-\s} (MV)^\ve (|t|+2)^\ve.  
\ena  
It follows that 
\bea\label{Sum/iii} 
- \int_{C_1} + \int_{C_2} + \int_{a-iV}^{a+iV}  
\ll (M^{1/4-\s}+MV^{-1}) (MV)^\ve (|t|+2)^\ve. 
\eea 
This finishes the estimation of the three integrals. 

It remains to compute $R$, the summation of the residues. 
The residue of the simple pole at $w=0$ is 
\bea\label{P/zeta} 
\f{P(s)}{\zeta(2s)}.
\eea  
For $\rho\not=0$, the pole 
at $w=\rho-s$ is simple with residue 
\bea\label{res/sim} 
c_q\f{M^{\rho-s}} {\zeta(2\rho)(\rho-s)}
\eea 
where $c_q$ is a constant depending on $q$. Before analyzing \eqref{res/sim} , we need a rough bound 
\bea\label{bd1/z}
\f{1}{|\zeta(2\rho)|}\ll_q 1. 
\eea 
This is trivially true if $|\rho|$ is bounded, say $|\rho|\leq 2$. For $|\rho|>2$ we apply the functional equation
$\zeta(s)=\g(s)\zeta(1-s)$ and that 
\bna 
|\g(s)| \sim \bigg(\f{|t|}{\pi}\bigg)^{1/2-\s} 
\ena 
as $|t|\to \infty$ in any strip $c_2\leq \s\leq c_3$ of finite width. 
Recall that all $\rho$ are lying on the imaginary axis and of the 
form of \eqref{P/pol}, so that 
\bna 
\f{1} {|\zeta(2\rho)|}\ll \f{1} {|\g(2\rho)| |\zeta(1-2\rho)|} \ll |\rho|^{\ve-1/2},  
\ena 
where we have also applied Lemma \ref{lem20}. This proves the claim \eqref{bd1/z}. 
Hence the summation of all the residues in \eqref{res/sim} with $|\Im\rho-t|\leq V$ is 
\bna
\sum_{|\Im\rho-t|\leq V} c_q \f{M^{\rho-s}} {\zeta(2\rho)(\rho-s)}
\ll M^{-\s} (|t|+V+2)^\ve \sum_{|\Im\rho-t|\leq V} \f{1} {|\rho-s|}. 
\ena  
To estimate the last sum over $\rho$, we note that, for $p|q$, 
\bna 
\#\bigg\{j: \f{H}{2}<\bigg|\f{2\pi j}{\log p}-t\bigg|\leq H\bigg\}\leq H\log p. 
\ena 
It follows that
\bna 
\sum_{|\Im\rho-t|\leq V} \f{1} {|\rho-s|}
&\leq& 1+\sum_{1\leq |\Im\rho-t|\leq V} \f{1} {|\Im\rho-t|} \\ 
&\ll& \sum_{p|q} \sum_{1\leq m\leq \log V/\log 2} \f{2^{m+1}\log p}{2^m}\\  
&\ll& (\log V)\sum_{p|q} (\log p)\\
&\ll&  (\log V)(\log q), 
\ena 
and hence  
\bea\label{S/res/sim} 
\sum_{|\Im\rho-t|\leq V} c_q \f{M^{\rho-s}} {\zeta(2\rho)(\rho-s)}
\ll M^{-\s} (|t|+V+2)^\ve
\eea  
on recalling that $q$ is a constant. The pole at $w=-s$ is of order $k$ where $k=\o(q)$, and therefore the residue is 
\bea\label{Mul/res} 
&& \lim_{w\to -s} \f{1}{(k-1)!} \f{d^{k-1}}{dw^{k-1}}\bigg((w+s)^k \f{P(s+w)} {\zeta(2s+2w)}\f{M^w}{w}\bigg) \nonumber\\  
&& \qquad \ll_q 
\sum_{0\leq l+m+n\leq k} \bigg| c_m 
\f{M^{-s}\log^l M }{s^{n}} \bigg| 
\ll M^{-\s+\ve},   
\eea 
where 
\bna 
c_m = \f{d^{m}}{ds^{m}} \bigg(\f{1}{\zeta(s)}\bigg)\bigg|_{s=0}.  
\ena 
It follows from \eqref{Mul/res}, \eqref{S/res/sim}, and \eqref{P/zeta} that 
\bea\label{R-ll} 
R-\f{P(s)}{\zeta(2s)} \ll M^{-\s+\ve}(|t|+V+2)^\ve, 
\eea 
and this finishes the treatment of $R$. 

Inserting \eqref{R-ll} and \eqref{Sum/iii} back to \eqref{V/iii} and \eqref{Mf/int/V}, we get 
\bna
\f{P(s)}{\zeta(2s)} -\sum_{m\leq M}^\infty \f{h(m)}{m^{s}}
\ll \f{M^{1+\ve}}{V} + \{ M V^{-1}  + M^{1/4-\s} + M^{-\s} \}(MV)^\ve (|t|+V+2)^\ve 
\ena 
which is, on taking $V=M^2$,  
\bna
\ll M^{1/4-\s} M^\ve (|t|+2)^\ve. 
\ena 
This proves the lemma for $M$ being half an odd integer.   
For general $M$ the above formula should be corrected by an error term of order 
$O(M^{-\s+\ve})$ by \cite[Lemma 3.19]{Tit}, and this correction can be absorbed. 
This completes the proof. 
\end{proof} 

\section{Proof of Theorem \ref{thm1}(ii)}
\setcounter{equation}{0}

Now we are ready to prove the conditional part of Theorem \ref{thm1}, 
applying the idea of Montgomery and Vaughan \cite{MonVau}. 
 
\begin{proof}[Proof of Theorem \ref{thm1}(ii)]  
Similarly to the unconditional case, we start from the formula 
\bna
M_f(x)
&=&\sum_{n\leq x} f(n) 
= \sum_{d m \leq x} \chi(d) h(m) \nonumber\\ 
&=& \bigg(\sum_{d m \leq x\atop m\leq M} + \sum_{d m\leq x\atop m>M}\bigg)\chi(d) h(m)\nonumber\\
&=& S_1 + S_2,  
\ena 
say. We assume in addition that $M<x$. For $S_1$, we invoke the trivial bound $|\sum_{d\leq x}\chi(d)|\leq q$ as well as  Lemma~\ref{lem/|h|}, to get 
\bea\label{S1/=} 
S_1
= \sum_{m\leq M} h(m) \sum_{d\leq x/m} \chi(d)   
\ll M^{1/2} \log^c x, 
\eea 
where $c=\o(q)$ and $q$ is treated as a constant. 

For $S_2$, we are going to apply analytic methods to prove in the following that 
\bea\label{S2/<<} 
S_2 
\ll M^{-1/4} x^{1/2+\ve}.   
\eea 
Once \eqref{S2/<<} is established, one sees immediately from this and \eqref{S1/=} 
that the optimal choice of $M$ is $M=x^{2/3}$, and this 
proves the Theorem \ref{thm1}(ii).  
 
\medskip 

It remains to show \eqref{S2/<<}.  
Again, we apply Perron's formula in the form of Lemma~\ref{lem1}, 
with $\psi(n)=n^\ve$ and $b= 1 +\f{1}{\log x}$. Then \eqref{sum|an} becomes, for $x$ being half an odd integer,  
\bea\label{Mfx=int} 
S_2  
= \f{1}{2\pi i} \int_{b-iT}^{b+iT} L(s, \chi) Z_M(s) \f{x^s}{s} ds +O\bigg(\f{x^{1+\ve}}{T}\bigg). 
\eea 
We shift the contour in \eqref{Mfx=int} to the left till the vertical line $\s=a$ where 
\bea\label{a=} 
a=\f12+\d 
\eea 
and $\d>0$ is arbitrarily small, without passing any singularities of the integrant in \eqref{Mfx=int}. 
The integral in \eqref{Mfx=int} can be written as 
\bea\label{i=iii} 
\int_{b-iT}^{b+iT} = - \int_{C_1} + \int_{C_2} + \int_{a-iT}^{a+iT}  
\eea 
where $C_1$ means the segment $\{s=\s-iT: a\leq \s\leq b\}$,   
and $C_2$ means $\{s=\s+iT: a\leq \s\leq b\}$. 
The integrals on $C_1$ and $C_2$ 
can be estimated to be small. Since we are not assuming the Riemann Hypothesis for 
$L(s,\chi)$, we invoke the convexity bound 
\bna 
L(s,\chi)\ll \{q(|t|+2)\}^{1/4}
\ena  
in the half-plane $\s\geq a$, and we can put $q$ into the implied constant.  
From this and Proposition \ref{prop/ZM}, we arrive at that   
\bna 
&&\bigg(-\int_{C_1} + \int_{C_2}\bigg) L(s,\chi) Z_M(s) \f{x^s}{s} ds \\
&& \ll \int_a^b M^{1/4} \bigg(\f{x}{M}\bigg)^\s (MT)^{\ve} T^{-3/4} d\s. 
\ena 
Note that $M<x$, and therefore the last quantity is 
\bna 
\ll (MT)^{-3/4} x (MT)^{\ve}.
\ena 
On the vertical segment, we apply Lemma~\ref{lem21}, so that  
\bna 
\int_{a-iT}^{a+iT} L(s,\chi) Z_M(s) \f{x^s}{s} ds 
&\ll& x^a M^{1/4-a+\ve}     
\int_{-T}^T |L(a+it,\chi)|  \f{(|t|+2)^\ve}{\sqrt{a^2+t^2}} dt \\ 
&\ll& M^{-1/4} x^{1/2}(MT)^\ve.  
\ena  
Inserting these back to \eqref{i=iii} and \eqref{Mfx=int}, we get 
\bna
S_2 
\ll \bigg(M^{-1/4} x^{1/2} + (MT)^{-3/4} x  + \f{x}{T}\bigg)(xT)^\ve.  
\ena 
Taking $T=x$ proves \eqref{S2/<<} for $x$ being half an odd integer, from which 
\eqref{S2/<<} for general $x$ follows easily. This completes the proof of 
Theorem \ref{thm1}(ii). 
\end{proof} 
   
\section{Proof of Theorem~\ref{thm3}}
\setcounter{equation}{0}
 
For simplicity we write 
\bea\label{Fs=} 
F(s)=\f{L(s, \chi) P(s)}{\zeta(2s)}; 
\eea 
the proof of Theorem~\ref{thm3} depends on the location of poles of $F$. 
Obviously poles of $F(s)$ must be zeros of $\zeta(2s)$. But we cannot exclude the possibility that 
some zeros of $\zeta(2s)$ may be cancelled by zeros of $L(s,\chi)$. It is this possibility that produces 
extra subtlety to the proof. What we indeed need to know is whether there are any simple zeros 
of $\zeta(2s)$, on the vertical line $\s=\f14$, that are not 
cancelled by the zeros of $L(s,\chi)$. The next lemma says that there are infinitely many of them.

\begin{lemma}\label{lem3} 
On the segment 
\bea\label{seg14} 
\{s=1/4+it: 0\leq t\leq T\},  
\eea   
there are $\gg T\log T$ simple zeros of $\zeta(2s)$ that are not zeros of $L(s,\chi)$.  
These zeros are simple poles of $F(s)$.  
\end{lemma}

Note that Lemma~\ref{lem3} does not depend on any unproved hypothesis.  
 
\begin{proof}  
It is known that a positive proportion of zeros of $\zeta(s)$ are simple and on the critical line $\s=\f12$, and 
Levinson \cite{Lev} has proved that this proportion is more than $\f13$. 
Thus $\zeta(2s)$ has $\gg T\log T$ simple zeros on the segment \eqref{seg14}. 

Next we are going to show that the number of zeros of $L(s,\chi)$ on the segment \eqref{seg14} 
is of lower order of magnitude. Let $\psi\bmod q$ denote a generic Dirichlet character, and $N(\a, T, \psi)$ the number of zeros of the attached $L$-function 
$L(s,\psi)$ lying in the rectangle 
\bna 
\a\leq \s\leq 1, \quad 0\leq t\leq T.
\ena 
Then Huxley's classical zero-density estimate \cite{Hux} states that  
\bna 
\sum_{\psi\bmod q} N(\a, T, \psi)\ll (qT)^{\f{12}{5}(1-\a)+\ve}, 
\ena 
where the summation is over all the characters $\psi$ modulo $q$. 
Now we take $\a=\f34$, and pick up just the term with $\psi=\chi$ on the left-hand side, 
to get 
\bna 
N(3/4, T, \chi)\ll (qT)^{3/5+\ve}. 
\ena 
In particular on the segment $\{s=\f34+it: 0\leq t\leq T\}$ the number of zeros of $L(s,\chi)$ 
is $\ll (qT)^{3/5+\ve}$. Note that our $\chi$ is real. By symmetry, on the segment \eqref{seg14},  
the number of zeros of $L(s,\chi)$ is also $\ll (qT)^{3/5+\ve}$. 

Hence after possible cancellations, there are still $\gg T\log T$ simple zeros of $\zeta(2s)$ left 
on the same segment 
\eqref{seg14}. This proves the lemma.    
\end{proof} 

\begin{proof}[Proof of Theorem~\ref{thm3}] 
The proof is naturally divided into two cases. 

(i) The first case is that there exists a constant $\th>\f14$ such that 
\bna 
\sup \{\b: F(s) \text{ has a pole at } \b+i\gamma\} = \th; 
\ena 
this implies that the Riemann Hypothesis for $\zeta$ is false. 
Suppose $M_f(x)\ll x^a$ for some $a$ satisfying $\f14<a<\th$. For $\s>1$ we apply partial integration to get  
\bea\label{1|z=} 
F(s)=\sum_{n=1}^\infty \f{f(n)}{n^s} 
= s \int_1^\infty \f{M_f(s)}{x^{s+1}} dx. 
\eea 
The last integral is absolutely convergent in the half-plane $\s\geq \f12 (a+\th)$, and hence represents 
an analytic function in the same region. This is a contradiction since $\f12 (a+\th)<\th$. This proves that  
\bna 
M_f(x)=\O(x^{a})
\ena 
for any $a$ satisfying $\f14<a<\th$, which is actually stronger than the statement of the theorem.

(ii) It remains to treat the case that $F(s)$ has no pole in the half-plane 
$\s>\f14$. Note that this does not mean the Riemann Hypothesis holds for $\zeta(2s)$, since its zeros 
off the line $\s=\f14$ may possibly be cancelled by the zeros of $L(s,\chi)$. 
In the present case \eqref{1|z=} holds in the half-plane $\s>\f14$. Thus 
\bea\label{45|con}
F(s)
= s \int_1^\infty \f{M_f(s)}{x^{s+1}} dx, 
\eea 
and finally we will make $\s\to\f14$ from the right.  
Suppose that 
\bna 
|M_f(x)|\leq \left\{
\begin{array}{lll}  
M_0, \quad &\text{for } 1\leq x\leq x_0, \\ 
\d x^{1/4},  \quad &\text{for } x\geq x_0.  
\end{array}
\right. 
\ena 
Then \eqref{45|con} gives 
\bna 
F(s)
&\leq& |s| M_0 \int_1^{x_0} \f{dx}{x^{5/4}} 
+ |s| \d  \int_{x_0}^\infty \f{dx}{x^{\s+3/4}} \\ 
&\leq& |s| M_0 \int_1^\infty \f{dx}{x^{5/4}} 
+ |s| \d  \int_1^\infty \f{dx}{x^{\s+3/4}} \\ 
&=& 4|s| M_0 +\f{|s|\d}{\s-\f14}. 
\ena 

Suppose that $\rho=\f14+i\gamma$ is a simple zero of $\zeta(2s)$ that is not cancelled 
by the zeros of $L(s,\chi)$; the existence of such a $\rho$ 
is guaranteed by Lemma~\ref{lem3}. This $\rho$ is a simple pole of $F(s)$, 
and we have $L(\rho,\chi)\not=0$. Of course $P(\rho)\not=0$. 
Putting $s=\s+ i\gamma$, and making $\s\to\f14$ horizontally from the right,   
\bna 
F(s) 
= \f{L(s,\chi)P(s)}{\zeta(2s)} 
= \f{L(s,\chi)P(s)}{\zeta(2s)-\zeta(2\rho)} 
\sim 
\f{L(\rho,\chi)P(\rho)}{2(\s-\f14) \zeta'(2\rho)}.  
\ena 
We therefore obtain a contradiction if 
\bna 
\d < \bigg|\f{L(\rho,\chi)P(\rho)}{4\rho \zeta'(2\rho)}\bigg|.  
\ena 
This proves the theorem. 
\end{proof} 

\medskip 

\noindent 
{\bf Acknowledgments.} The author would like to record her sincere thanks to  
Liming Ge, Fuzhou Gong, and Nanhua Xi for supports and encouragements, to Marco 
Aymone and Bingrong Huang for helpful suggestions leading to improvements in Theorem \ref{thm1}(ii),    
and to Oleksiy Klurman for bibliographical information.

\end{document}